\documentclass[12pt]{article}
\usepackage{amsmath,amssymb,amsfonts,graphicx}
\makeatletter \@addtoreset{equation}{section} \makeatother
\newtheorem{thm}{Theorem}[section]

\newtheorem{rem}{Remark}[section]

\begin{document}
\begin{center}

{\bf \large On Products of Random Matrices and certain Hecke
Algebras associated with Groups of $2\times 2$ Matrices}

\vspace*{1cm} {\bf Jafar Shaffaf}\footnote{Institute for Studies
in Theoretical Physics and Mathematics (IPM) and Sharif University
of Technology, Tehran, Iran.  Email : shaffaf@ipm.ir}
\vspace{0.3cm}
\end{center}

\vspace*{1cm}

\begin{abstract}
\noindent  The determination of the density functions for products
of random elements from specified classes of matrices is a basic
problem in random matrix theory and is also of interest in
theoretical physics. For connected simple Lie groups of $2\times
2$ matrices and conjugacy and spherical classes a complete
solution is given here. The problem/solution can be re-stated in
terms of the structure of certain Hecke algebras attached to
groups of $2\times 2$ matrices.
\end{abstract}

\begin{center}
\end{center}

\section{Introduction}

Let $(G,K)$ be a symmetric pair where $G$ is a semi-simple Lie
group with finite center.  For $G$ compact let ${\mathcal
H}_c={\mathcal H}_c(G)$ be the algebra under convolution generated
by the invariant measures concentrated on conjugacy classes in
$G$, and ${\mathcal H}_s={\mathcal H}_s(G,K)$ be the convolution
algebra generated by the invariant measures on spherical classes
${\mathcal O}_a=KaK\subset G$.  It is elementary that ${\mathcal
H}_c$ and ${\mathcal H}_s$ are commutative algebras with unit.
These algebras arise naturally in random matrix theory, and
theoretical physics. In fact, the product structures for
generators of these algebras are given by the density functions of
products of random matrices chosen according to invariant measures
on conjugacy or spherical classes; and in string theory generic
D-branes are localized along the product of (twisted) conjugacy
classes of the Lie group.  We will not discuss specific physical
interpretations in this paper (except for a remark in \S 5), and
the reader is referred to [Q], and references thereof, for a
detailed discussion of the issues of interest in string theory.

This work may be regarded as the initial attempt at understanding
the structure of the Hecke algebras ${\mathcal H}_c$ and
${\mathcal H}_s$ and considers only the special case of $2\times
2$ matrices. The support of the density function for products of
two generators of ${\mathcal H}_c(SU(n))$ is determined in [AW]
and is exhibited by a set of linear inequalities on the Lie
algebra of a maximal torus. The result in [AW] is essentially a
reformulation of a theorem about singular connections on a
holomorphic vector bundle on a Riemann surface with marked points
[B], or equivalently a theorem of Mehta and Seshadri in algebraic
geometry.  Neither theorem is applicable for the computation of
the density function for the product of two generators of the
Hecke algebra.

A complete structure theorem for ${\mathcal H}_c(SU(2))$ is given
in Theorem \ref{thm:shaffaf1} below.  For spherical classes
attached to groups of $2\times 2$ matrices we consider the
symmetric pairs $(SU(2),S(U(1)\times U(1))$, $(SL(2,\mathbb{R}),
SO(2))$ and $(SL(2,\mathbb{C}),SU(2))$.  Theorem
\ref{thm:shaffaf2} gives density functions for products of
spherical classes.  In the final section we present some numerical
results in relation to products of conjugacy classes in $SU(2)$.
Each conjugacy class, with the induced metric, is a copy of $S^2$
equipped with a metric of constant curvature.  It is noticed
(perhaps surprisingly) that if $S^2$ is dicretized according to
the prescription of Thomson's problem (minimizing Coulomb
potential or intuitively "the best equally spaced distribution")
then convergence to the predicted measure is much slower than if
the points were chosen randomly.

In [JM] the question of whether a product of conjugacy classes
${\mathcal C}_{\alpha_1},\ldots,{\mathcal C}_{\alpha_n}$ contains
the identity element and/or is the entire group $SU(2)$ is
studied. By successive applications of Theorem \ref{thm:shaffaf1}
of this paper one can recover the results in [JM] and in fact give
a more precise version of it.  This subject has not been
elaborated on here.  The argument in [JM] depends on the results
in [AW], but the proofs presented here are self-contained and of a
more elementary nature.

The author wishes to thank A. Katanforoush for the numerical
results, R. Szabo, M. Sheikh-Jabbari, S. Shahshahani and
especially Mehrdad Shahshahani for many stimulating discussions.

\section{Statement of main results}
A conjugacy class $C_\theta\subset SU(2)$ is uniquely determined
by the eigenvalues $e^{\pm i\theta}$ of a matrix in $C_\theta$.
Each conjugacy class is a homogeneous space, and for $C_\theta$
this measure $\mu_\theta$ is uniquely normalized to be $4\pi
\sin^2\theta$ in accordance with Weyl's integration formula and
the normalization ${\rm vol.}(SU(2))=4\pi^2$.

\begin{thm}
\label{thm:shaffaf1} Let $\mu_\alpha$ and $\mu_\beta$ be the
invariant measure on the conjugacy classes $C_\alpha$ and
$C_\beta$ respectively (regarded as singular distribution on $G$).
Then $\mu_\alpha\star\mu_\beta$ is an absolutely continuous
conjugation invariant measure on $G$ relative to the Haar measure
and its density is given by
\begin{equation}
\label{eq:shaffaf1} \mu_\alpha\star\mu_\beta =
\begin{cases}4\pi^2\sin\alpha\sin\beta\sin\theta, & {\rm for} ~~\alpha-\beta\le \theta \le
\alpha+\beta,\cr -4\pi^2\sin\alpha\sin\beta\sin\theta, & {\rm for}
~~-\alpha-\beta\le \theta \le -\alpha+\beta,\cr 0,& {\rm
otherwise},
\end{cases}
\end{equation}
where $0<\beta\le \alpha <\pi$.
\end{thm}

Let $K=S(U(1)\times U(1))\subset SU(2)$, then $(G,K)$ is a
symmetric pair of compact type [H1]. Orbits of the action of
$K\times K$ on $G$ via
\begin{eqnarray*}
a\longrightarrow kak^\prime,~~~~(k,k^\prime)\in K\times K,~~a\in G
\end{eqnarray*}
are called spherical classes and are denoted by ${\cal O}_a$.
${\mathcal O}_a$ is a homogeneous space and it is a simple
calculation that $|a_{11}|$ is constant on ${\cal O}_a$ and
uniquely determines it. Here $a=\begin{pmatrix} a_{11}&a_{12}\cr
a_{21}&a_{22}
\end{pmatrix}\in SU(2)$ is any matrix in ${\cal O}_a$.
Each spherical class ${\mathcal O}_a$ carries an invariant measure
$\lambda_a$.  The total mass of the measure $\lambda_a$ of a
spherical class ${\cal O}_a$ is $2|a_{11}|$, and we set
$r(a)=|a_{11}|$.

In the non-compact cases $(SL(2,\mathbb{R}), SO(2))$ and
$(SL(2,\mathbb{C}),SU(2))$ the real Cartan subgroup $A$ is
\begin{eqnarray*}
A = \{ a_t = \begin{pmatrix} e^{\frac{t}{2}} & 0 \\ 0 &
e^{-\frac{t}{2}}
\end{pmatrix}~|~ t\in \mathbb{R}\}
\end{eqnarray*}
$A_+=\{a_t\in A~|~t>0\}$ parametrizes spherical classes. The Haar
measure for the Cartan (polar) decomposition $G \simeq KAK$ is
\begin{equation}
\label{eq:Cartan1} \int_G f(g) dg = c \int_K \int_K
\int_0^{\infty} f(k_1 a k_2 ) \delta(t) dt  d k_1  dk_2
\end{equation}
where $c$ is a suitable constant (see [H2]), $\delta(t)=\sinh^2
t$, and $\epsilon=1~{\rm or}~2$ according as $G=SL(2,\mathbb{R})$
or $SL(2,\mathbb{C})$. The volume (area) of the spherical class
${\mathcal O}_{a_t}$ is
\begin{equation}
\label{eq:shaffaf22} {\rm vol.}(\mathcal{O}_{a_t}) = [4 \pi^2 c
\sinh t]^\epsilon.
\end{equation}

\begin{thm}
\label{thm:shaffaf2} The product formula for spherical classes in
the three cases of symmetric pairs attached to groups of $2\times
2$ matrices are:

\noindent {\rm (A)}- Let $\lambda_a$ and $\lambda_b$ be two
(singular) spherical measures concentrated on the spherical
classes ${\cal O}_a$ and ${\cal O}_b$ respectively. Then
$\lambda_a\star\lambda_b$ is absolutely continuous relative to the
Haar measure on $SU(2)$ and given by
\begin{equation}
\label{eq:Shaffaf2}
\begin{cases}\frac{16 \pi^2|a_{11}
b_{11}|u}{\sqrt{c_1^2-(u^2-c_0)^2}},&{\rm for} ~~\sqrt{c_0-c_1}\le
u\le \sqrt{c_0+c_1},\cr 0,&{\rm otherwise},
\end{cases}
\end{equation}
where $c_0$ and $c_1$ are symmetric functions of $a$ and $b$ and
given by
\begin{eqnarray*}
c_0 &=& r^2(a) r^2(b)+ (1- r^2(a))(1 - r^2(b)),\\  c_1 &=& 2r(a)
r(b) \sqrt{(1- r^2(a))(1 - r^2(b))}.
\end{eqnarray*}

\noindent {\rm (B)} - Let $\lambda_{a_{t_1}}$ and
$\lambda_{a_{t_2}}$ be the (singular) invariant measures
concentrated on the spherical classes ${\cal O}_{a_{t_1}}$ and
${\cal O}_{a_{t_2}}$ in $SL(2,\mathbb{R})$. Then
$\lambda_{a_{t_1}}\star\lambda_{a_{t_2}}$ is spherical and
absolutely continuous relative to the Haar measure, and for a
continuous spherical function $f$ on $SL(2, \mathbb{R})$ we have
\begin{eqnarray*}
\lambda_a\star\lambda_b(f) = 4 c^2 \pi^2 \sinh t_1 \sinh t_2 \int_
{I_{t_1, t_2}} f(r)  \frac{ \sinh r}{\sqrt{ c_2^2 - ( c_1 - \cosh
r)^2}} ~ dr ~,
\end{eqnarray*}
where $ I_{t_1, t_2} = [t_2 - t_1, t_2 + t_1 ]$, $c_1 = \cosh t_1
\cosh t_2$, and $c_2 = \sinh t_1 \sinh t_2 $.

\noindent {\rm (C)} - Let $\lambda_{a_{t_1}}$ and
$\lambda_{a_{t_2}}$ be the (singular) invariant measures
concentrated on the spherical classes ${\cal O}_{a_{t_1}}$ and
${\cal O}_{a_{t_2}}$ in $SL(2,\mathbb{C})$. Then
$\lambda_{a_{t_1}}\star\lambda_{a_{t_2}}$ is spherical and
absolutely continuous relative to the Haar measure, and for a
continuous spherical function $f$ on $SL(2, \mathbb{C})$ we have
\begin{eqnarray*}
\lambda_{a_{t_1}}\star\lambda_{a_{t_2}}(f) = 32 c^2 \pi^6 \sinh
t_1\sinh t_2 \int_ {I_{t_1, t_2}} f(r) \sinh r dr ~,
\end{eqnarray*}
where $ I_{t_1, t_2} = [t_2 - t_1, t_2 + t_1 ]$.
\end{thm}

It may be of interest to normalize the measures on the spherical
classes to probability measures and determine the corresponding
empirical measure of products.  For such a normalization the
density functions determined in Theorem \ref{thm:shaffaf2} become
\begin{equation}
\label{eq:shaffaf81} {\rm A:}~~\frac{1}{2\pi} \frac {u}{
\sqrt{c_1^2-(u^2 - c_0)^2}},~~~{\rm B:}~~\frac{1}{\pi} \frac{
\sinh r}{\sqrt{ c_2^2 - ( c_1 - \cosh r)^2}},~~~{\rm
C:}~~\frac{\sinh r}{2 \sinh t_1 \sinh t_2}.
\end{equation}

\begin{rem}
\label{rem:Shaffaf1}{\rm Let $\tilde{{\mathcal H}_c}$ and
$\tilde{{\mathcal H}_s}$ denote the completions of ${\mathcal
H}_c$ and ${\mathcal H}_s$ in the weak topology.  Then in all
cases considered $\tilde{{\mathcal H}_c}$ and $\tilde{{\mathcal
H}_s}$ contain the corresponding $L^1$ space as a dense ideal.
Furthermore, by the above analysis, for every pair of generators
$\mu_1$ and $\mu_2$, we have $\mu_1\star\mu_2\in L^1$.}
\end{rem}

\section{Proof of Theorem \ref{thm:shaffaf1}}

First we show that $\mu_\alpha\star\mu_\beta$ is an $L^p$ function
for $p\le 2$.  The Fourier expansion of the singular measure
$\mu_\alpha$ is given by

\begin{equation}
\label{eq:shaffaf3} \sum_{\rho\in \hat{G}} d_\rho {\rm
Tr}\big[\rho (\mu_\alpha)\rho (g)\big],
\end{equation}
where $\hat{G}$ is the set of irreducible representations of $G$
and the Fourier transform of a measure $\mu$ is defined as
\begin{eqnarray*}
\rho(\mu)=\int \rho(g^{-1})d\mu.
\end{eqnarray*}
This series does not converge in the ordinary sense of convergence
of series of functions since the measure $\mu_\alpha$ is singular,
however, it converges in the weak sense.  The convolution product
$\mu_\alpha\star\mu_\beta$ can be calculated from
\begin{equation}
\label{eq:shaffaf4} \sum_{\rho\in \hat{G}} d_\rho {\rm
Tr}\big[\rho (\mu_\alpha)\rho(\mu_\beta)\rho (g)\big].
\end{equation}
Since the character $\chi_\rho(g)$ is independent of $g\in
C_\alpha$, and is given by $\frac{\sin (n+1)\alpha}{\sin \alpha}$
for $\rho$ is the symmetric $n^{{\rm th}}$ representation of
$G=SU(2)$,
\begin{eqnarray*}
\rho(\mu_\alpha) = \frac{4\sin^2\alpha}{d_\rho}{\rm Tr}(\rho(g))I=
\frac{4\chi_\rho(C_\alpha)\sin^2\alpha }{d_\rho}I.
\end{eqnarray*}
Applying the Plancherel theorem to $\mu_\alpha\star \mu_\beta$, we
obtain
\begin{eqnarray*}
||\mu_\alpha\star \mu_\beta||^2 &=&\sum \frac{16
\sin^2\alpha\sin^2\beta
\chi_\rho(C_\alpha)\chi_\rho (C_\beta)}{d_\rho^4} d^2_\rho \\
&= & \sum_{n\ge 1} \frac{16}{n^2}
\sin^2(n+1)\alpha\sin^2(n+1)\beta.
\end{eqnarray*}
This series converges absolutely and therefore
$\mu_\alpha\star\mu_\beta$ is a square integrable function.  From
the Cauchy-Schwartz inequality it follows that
$\mu_\alpha\star\mu_\beta$ is absolutely integrable.

Therefore the measure $\mu_\alpha\star\mu_\beta$ is a conjugation
invariant function and it can be interpreted as the (defective)
density function for the space of solutions $c$ to the equation
\begin{eqnarray*}
abc=e,~~~{\rm where}~~a\in C_\alpha,~b\in C_\beta.
\end{eqnarray*}
Since this density is conjugation invariant we represent it as the
function $\nu (\theta)$ on the maximal torus of diagonal matrices
given by
\begin{eqnarray*}
\nu(\theta)= ((\mu_\alpha\star\mu_\beta)\star \mu_\theta)(e).
\end{eqnarray*}
Let ${\cal R}$ denote the left regular representation of $G$, and
$\Delta_\alpha= 4\pi\sin^2\alpha$ be the factor appearing in Weyl
integration formula for conjugacy invariant functions.  The
Fourier transform of a function $\psi$ on $G$ at a representation
$\rho$ is defined as
\begin{eqnarray*}
\rho (\psi)= \int \rho(x^{-1}) \psi (x) dx.
\end{eqnarray*}
Let $d_\rho$ denote the dimension of the representation $\rho$,
$\chi_\rho$ its character and $\hat{G}$ the space of (complex)
irreducible representations of $G$. Taking Fourier transform and
decomposing the regular representation ${\cal R}$ of $G$ in the
usual way we obtain
\begin{eqnarray*}
\nu(\theta)&=& ((\mu_\alpha\star\mu_\beta)\star \mu_\theta)(e)\\
&=& \frac{1}{{\rm vol.}(G)} {\rm Tr} {\cal
R}((\mu_\alpha\star\mu_\beta)\star \mu_\theta)\\
&=& \frac{1}{{\rm vol.}(G)}~\sum_{\rho\in \hat{G}}d_\rho{\rm Tr} \rho((\mu_\alpha\star\mu_\beta)\star \mu_\theta)\\
&=& \frac{1}{{\rm vol.}(G)}~\sum_{\rho\in \hat{G}} d_\rho
\frac{\Delta_\alpha\Delta_\beta\Delta_\theta}{d_\rho^3}\chi_\rho(\alpha)\chi_\rho(\beta)\chi_\rho(\theta){\rm
Tr}(I)\\
&=&\frac{\Delta_\alpha\Delta_\beta\Delta_\theta}{{\rm vol.}(G)}~
\sum_{\rho\in
\hat{G}}\frac{\chi_\rho(\alpha)\chi_\rho(\beta)\chi_\rho(\theta)}{d_\rho}.
\end{eqnarray*}
Therefore
\begin{equation}
\label{eq:shaffaf6} \nu(\theta)= 16 \pi
\sin^2\alpha\sin^2\beta\sin^2\theta~ \sum_{\rho\in
\hat{G}}\frac{\chi_\rho(\alpha)\chi_\rho(\beta)\chi_\rho(\theta)}{d_\rho}
\end{equation}

Irreducible representations of $G=SU(2)$ are determined by
dimension $k\ge 2$. The character of the $k$-dimensional
representation $\rho_k$ is:
\begin{eqnarray*}
\chi_k (g)=\sum_{j=0}^{k}
e^{(j-2)i\theta}=\frac{\sin(k+1)\theta}{\sin\theta}
\end{eqnarray*}
where $e^{\pm i\theta}$ is the eigenvalue of the conjugacy class
of $g$.  Substituting in (\ref{eq:shaffaf6}) we obtain
\begin{equation}
\label{eq:shaffaf7} \nu(\theta) = 16 \pi
\sin\alpha\sin\beta\sin\theta~ \sum_{k=0}^\infty \frac{\sin (k+1)
\alpha  \ \sin (k+1) \beta  \ \sin (k+1)\theta  } {k+1}.
\end{equation}
Now let $f(\theta)$ denote the function defined by
(\ref{eq:shaffaf1}).  Since $f$ is an even function its Fourier
expansion is of the form $f(\theta)=\frac{a_0}{2}+\sum a_n\cos
nx$, and
\begin{eqnarray*}
a_n=\frac{2}{\pi}\big[\frac{1}{n+1}\sin(n+1)\alpha\sin(n+1)\beta -
\frac{1}{n-1}\sin(n-1)\alpha\sin(n-1)\beta\big]
\end{eqnarray*}
Substituting in the Fourier expansion we obtain
\begin{eqnarray*}
f(\theta)= 8\pi \sin\alpha
\sin\beta\sum_{m=0}(\frac{\sin(m+1)\alpha\sin(m+1)\beta}{m+1}-
\frac{\sin(m-1)\alpha \sin(m-1)\beta}{m-1})\cos m\theta.
\end{eqnarray*}
Using the elementary identity $\cos m\theta-\cos
 (m+2)\theta  = 2\sin(m+1)\theta\sin\theta$
the series becomes telescopic and simplifies to
\begin{eqnarray*}
f(\theta)= 16\pi \sin\alpha
\sin\beta\sin\theta\sum_{m=0}\frac{\sin(m+1)\alpha\sin(m+1)\beta\sin(m+1)\theta}{m+1},
\end{eqnarray*}
which is identical with (\ref{eq:shaffaf7}). $\blacksquare$

\begin{rem}
\label{rem:shaffaf1} {\rm In the above proof we made use of the
fact that the measure $\mu_\alpha\star\mu_\beta$ is a conjugation
invariant function that can be interpreted as the (defective)
density function for the space of solutions $c$ to the equation
\begin{eqnarray*}
abc=e,~~~{\rm where}~~a\in C_\alpha,~b\in C_\beta
\end{eqnarray*}
and this density is given by (\ref{eq:shaffaf6}).  For the case of
finite groups this formula is well-known ([S], p.68). }
\end{rem}

\begin{rem}
\label{rem:shaffaf2} {\rm It is possible to prove Theorem
\ref{thm:shaffaf1} without the use of harmonic analysis and by
integral formulae similar to those used for the proof of Theorem
\ref{thm:shaffaf2} below. However, it appears that the above
argument is possibly generalizable to products of conjugacy
classes in compact connected semi-simple Lie groups, but the one
based on integral formulae is not. }
\end{rem}

\section{Proof of Theorem \ref{thm:shaffaf2}}

Since the proofs of (B) and (C) are essentially the same
computation we prove (A) and (C) only.  Introduce coordinates on
$SU(2)$ by:
\begin{equation}
\label{eq:shaffaf24} (\rho, \varphi,
\psi)\longrightarrow\begin{pmatrix} \rho e^{i\varphi}
 & \sqrt{1-\rho^{2}}e^{-i\psi}
\\-\sqrt{1-\rho^{2}}e^{i\psi}  &
\rho e^{-i\varphi}\end{pmatrix}, ~~~(\rho,\varphi,\psi)\in
[0,1]\times [0,2\pi]\times [0,2\pi]
\end{equation}
The Haar measure on $SU(2)$ in the $(\rho, \varphi, \psi)$-
coordinates is easily calculated by computing $g^{-1}dg$, a basis
of left invariant 1-forms $\omega_1,\omega_2,\omega_3$ and then
taking their wedge product to obtain:
\begin{eqnarray*}
\Omega= \omega_1\wedge\omega_2\wedge\omega_3= 2\rho  d \rho d
\varphi  d \psi.
\end{eqnarray*}
With this normalization ${\rm vol}(SU(2))=4\pi^2$ as before.

Since both $f$ and $\lambda_a$ are $K$-bi-invariant,
$\lambda_a\star f(x)$ is $K$-bi-invariant and therefore to compute
$\mu_a\star f(x)$ we can assume that $x$ is of the form
\begin{equation}
\label{eq:shaffaf11} x= \begin{pmatrix} t & \overline{w}
\\ -w & t\end{pmatrix}
\end{equation}
where $t$ is a real number and $w=s e^{i\alpha}$ a complex number
with $t^2+s^2=1$. In $(\rho, \varphi, \psi)$- coordinates on
$SU(2)$ we have

\begin{eqnarray*}
\lambda_a\star\check{f}(x)= \int_{{\cal O}_a} f(y x^{-1}) dy=
2r(a) \int_0^{2\pi}\int_0^{2\pi} f( y(\rho, \varphi, \psi)x^{-1})
 d\varphi d\psi
\end{eqnarray*}
With $x$ represented as in (\ref{eq:shaffaf11}) we have
\begin{eqnarray*}
\lambda_a \star \check {f}(x)&= & 2r(a) \int_0^{2\pi}\int_0^{2\pi}
f \left ( \begin{pmatrix}
 \rho e^{i\varphi}
 & \sqrt{1-\rho^{2}}e^{-i\psi}
\\
-\sqrt{1-\rho^{2}}e^{i\psi}  & \rho e^{-i\varphi}
\end{pmatrix}\begin{pmatrix} t & \overline{w}
\\ -w & t  \end{pmatrix} \right ) d\varphi
d\psi \\
 &=& 2r(a) \int _0^{2\pi} \int_
0^{2\pi} f \left (\begin{pmatrix} r(a)t e^{i \varphi} + b
\sqrt{1-\rho^2} e^{-i\psi}
 & \star \\
\star  & \star \end{pmatrix} \right ) d \varphi d\psi
\end{eqnarray*}
Since $f$ is spherical, it depends only on the norm of the $(1,1)$
entry of the above matrix. The square of the norm of the $(1,1)$
entry is

\begin{eqnarray*}
t^2 \ r^2(a) + s^2 (1 - r^2(a))+ 2t  s r(a) \sqrt{1- r^2(a)}
\cos(\varphi + \psi- \alpha).
\end{eqnarray*}
Substituting $ t=r(x)$ and $s=\sqrt{1-r^2(x)}$, the norm of the
$(1,1)$ entry becomes
\begin{eqnarray*}
 |r(a) te^{i\varphi} + b
\sqrt{1- \rho^2} e^{-i\psi}|= \sqrt{c_0 + c_1\cos(\varphi + \psi-
\alpha)}.
\end{eqnarray*}
Therefore
\begin{eqnarray*}
 \lambda_a\star \check{f}(x)=2r(a)
\int _0^{2\pi} \int_ 0^{2\pi} f(\sqrt{c_0 + c_1\cos(\varphi +
\psi- \alpha)}) \ d\varphi \ d\psi.
\end{eqnarray*}
The change of variable
\begin{eqnarray*}
(u,v)= (\sqrt{c_0 + c_1\cos(\varphi + \psi- \alpha)}, \psi),
\end{eqnarray*}
is a 2 to 1 covering.  Its Jacobian is given by
\begin{eqnarray*}
\frac{\partial(\varphi,\psi)}{\partial(u, v)}=\frac{2u}{c_1
\sin(\varphi + \psi- \alpha)} =\frac{2u}{\sqrt{c_1^2 - (
u^2-c_0)^2}}
\end{eqnarray*}
Therefore
\begin{eqnarray*}
\lambda_a \star \check {f}(x) &=& 4 \pi r(a)\int _0^{2\pi} \int_
{\sqrt{c_0 - c_1}}^{\sqrt{c_0 + c_1}} f(u) \frac {u}{
\sqrt{c_1^2-(u^2 - c_0)^2}} \  du dv \\
&=& 2r(a) \int_ {\sqrt{c_0 -c_1}}^{\sqrt{c_0+ c_1}} f(u) \frac
{u}{ \sqrt{c_1^2-(u^2 - c_0)^2}} \  du.
\end{eqnarray*}
To compute $\lambda_a \star \lambda_b(f)$ we set $g(x) = \mu_b
\star \check{f}(x)$. Then $g$ is spherical and
\begin{eqnarray*}
\lambda_a\star\lambda_b(f) & =&\lambda_a \star
(\lambda_b \star\check{f})(e)\\&=&(\lambda_a \star g)(e)\\
&=&\int_{{\cal O}_a} g(x)d\lambda_a (x)\\&=& g(a)
\textrm{vol.}({\cal O}_a),
\end{eqnarray*}
Now
\begin{eqnarray*}
{\rm vol.}({\cal O}_a)=2r(a)\int_0^{2\pi}\int_0^{2\pi} \ d \varphi
\ d\psi = 8\pi^2 r(a).
\end{eqnarray*}
Therefore
\begin{eqnarray*}
\mu_a\star\mu_b(f)= 8\pi^2 r(a)  g(a).
\end{eqnarray*}
Substituting from the calculation of
$g(a)=\lambda_b\star\check{f}(a)$ above we obtain
\begin{eqnarray*}
\lambda_a\star\lambda_b(f)= g(a) \ \textrm{vol.}({\cal O}_a) =
16\pi^2 r(a) r(b) \int_{\sqrt{c_0 - c_1}}^{\sqrt{c_0 + c_1}} f(u)
\ \frac {u}{ \sqrt{c_1^2-(u^2 - c_0)^2}} \  du
\end{eqnarray*}
This completes the proof of part (A).

\noindent {\bf Proof of part (C)} - Since both $f$ and $\lambda_a$
are $K$-bi-invariant, $\lambda_a\star f(x)$ is $K$-bi-invariant
and therefore to compute $\mu_a\star f(x)$ we can assume that $x$
is of the form $x=\begin{pmatrix}
  e^{\frac{t}{2}} & 0 \\ 0 &
e^{-\frac{t}{2}}
\end{pmatrix}$.
As before let $\{\theta_n\}$ be a sequence of spherical functions
on $G$ converging weakly to the (singular) invariant measure
$\lambda_a$ on the orbit $\mathcal{O}_a$. Applying the polar
coordinate decomposition, for the convolution $\lambda_a \star
\check{f}(x)$ we have
\begin{eqnarray*}
\lambda_a \star \check{f}(x) &=& \int_{\mathcal{O}_a} f(y x^{-1}) dy
\\ &=& \lim_{n \rightarrow \infty} \int_G \theta_n(g) f(g x^{-1} ) dg
\\&=&c \lim_{n \rightarrow \infty} \int_K \int_K \int_A \theta_n(k_1 a' k_2 ) f ( k_1 a' k_2
x^{-1}) \delta(a') da' dk_1 dk_2 \\& =& c\lim_{n \rightarrow
\infty
} \int_K \int_A \theta_n(a') f(a' k x^{-1}) \delta(a') da' dk \\
&=& c \delta(t_1) \int_K f( a k x^{-1}) dk
\end{eqnarray*}
Writing $M = a k x^{-1} = k_1 a_1 k_2$ , where $k_1, k_2 \in K =
SU(2)$ and using the coordinates in (\ref{eq:shaffaf24})
\begin{eqnarray*}
 a_1 = \begin{pmatrix}e^{\frac{r}{2}} & 0 \\
0 & e^{\frac{-r}{2}}
\end{pmatrix} ,
~~ k = \begin{pmatrix} \rho e^{i\varphi}
 & \sqrt{1-\rho^{2}}e^{-i\psi}
\\-\sqrt{1-\rho^{2}}e^{i\psi}  &
\rho e^{-i\varphi}\end{pmatrix},
\end{eqnarray*}
we compute $r$ in term of $t_1$, $t_2$ and $k$:
\begin{equation}
\label{eq:shaffaf13} 2 \cosh r={\rm Tr}(a_1^2)={\rm Tr }( k_1
a_1^2 k_1^{-1})={\rm Tr}(M M^*).
\end{equation}
On the other hand we have
\begin{eqnarray*}
M = a k x^{-1} = \begin{pmatrix} e^{\frac{1}{2}(t_1 + t)}\rho
e^{i\varphi}
 &e^{ \frac{1}{2}(t_1 - t)}\sqrt{1-\rho^{2}}e^{-i\psi}
\\-e^{\frac{1}{2}(t - t_1)}\sqrt{1-\rho^{2}}e^{i\psi}  &
e^{-\frac{1}{2}(t_1 + t)}\rho e^{-i\varphi}\end{pmatrix}
\end{eqnarray*}
Therefore
\begin{eqnarray*}
{\rm Tr}(MM^*) = 2\rho^2 \cosh (t + t_1) + 2(1-\rho^2) \cosh (t -
t_1)~ ,
\end{eqnarray*}
Comparing with (\ref{eq:shaffaf13}) we obtain after a simple
calculation
\begin{eqnarray*}
\cosh r = \cosh t_1 \cosh t - (2 \rho^2 - 1)\sinh t_1 \sinh t
\end{eqnarray*}
Now set
\begin{eqnarray*}
c_1 = \cosh t_1 \cosh t ,~ c_2 = \sinh t_1 \sinh t
\end{eqnarray*}
The function $f$ is spherical so it only depends on the component
$r$ and therefore
\begin{eqnarray*}
\lambda_a \star \check{f}(x)= c \delta(t_1) \int_K f(a k x^{-1}) dk
= c \delta(t_1) \int_0^{2\pi}\int_0^{2\pi} \int_0^1 f(r) (2\rho)
 d\rho d\varphi d\psi
\end{eqnarray*}
Now make the change of coordinate $ \cosh r = c_1 - c_2 (2 \rho^2
-1)$ and note that, assuming that $t>t_1$, $r$ ranges over
$I_{t,t_1}=[t-t_1, t+t_1]$ as $\rho$ ranges over $[0,1]$.
Substituting in the above integral for $\lambda_a\star
\check{f}(x)$ we obtain
\begin{equation}
\label{eq:shaffaf26} \lambda_a \star \check{f}(x)= \frac{2 \pi^2
c}{c_2} \delta(t_1) \int_{I_{t,t_1}} f(r) \sinh r dr ~ ,
\end{equation}
To compute $\lambda_a \star \lambda_b(f)$ we set $g(x) = \mu_b
\star \check{f}(x)$. Then $g$ is spherical and
\begin{eqnarray*}
\lambda_a\star\lambda_b(f) & =&\lambda_a \star
(\lambda_b \star\check{f})(e)\\&=&(\lambda_a \star g)(e)\\
&=&\int_{{\cal O}_a} g(x)d\lambda_a (x)\\&=& g(a)
\textrm{vol.}({\cal O}_a),
\end{eqnarray*}
Using (\ref{eq:Cartan1}) one easily obtains
\begin{equation}
\label{eq:shaffaf23} {\rm vol.}({\mathcal O_{a_t}})=c({\rm vol.}
(SU(2)))^2\sinh^2t=16\pi^4\sinh^2t.
\end{equation}
Substituting from (\ref{eq:shaffaf23}) we obtain
\begin{eqnarray*}
\lambda_a\star\lambda_b(f)= 16 c \pi^4 (\sinh t_1)^2  g(a).
\end{eqnarray*}
Equation (\ref{eq:shaffaf26}) implies
\begin{eqnarray*}
\lambda_a\star\lambda_b(f) = 32 c^2 \pi^6  \sinh t_1 \sinh t_2 \int_
{I_{t_1, t_2}} f(r) \sinh r dr ~,
\end{eqnarray*}
which completes the proof of the theorem. $\blacksquare$

\section{Discretization and Numerical Results}

The discretization and quantization of conjugacy classes are
essentially different issues, and the latter is based on the
Kirillov orbit method (see [K]).  It is of interest in theoretical
physics and is also relevant to the subject matter of this paper.
The dual of the Lie algebra of $SU(2)$ is identified with the set
of $2\times 2$ skew hermitian matrices of trace 0:
\begin{eqnarray*}
\xi (a,b,c)=\begin{pmatrix}ic&a+ib\cr -a+ib&-ic  \end{pmatrix}
\end{eqnarray*}
Under conjugation action of $G$ the orbits are the spheres
$\Sigma_r=\{\xi (a,b,c)~|~a^2+b^2+c^2 = r^2\}$.  It is customary
in physics to assign the symmetric $n^{{\rm th}}$ power
representation $\rho_n$ of $SU(2)$ to the sphere $\Sigma_r$ of
area $n\in\mathbb{Z}_+$. According to the Clebsch-Gordon formula,
$\rho_n\otimes\rho_m$ decomposes as
\begin{eqnarray*}
\rho_n\otimes \rho_m\simeq \sum_{|n-m|}^{n+m} \rho_k,~~~{\rm
where}~~k\equiv m+n~~{\rm mod}~2.
\end{eqnarray*}
On the other hand the Minkowski sum of spheres of radii $r_1$ and
$r_2$ in $\mathbb{R}^3$ is precisely the spherical shell defined
by $|r_2-r_1|\le ||\xi||\le r_1+r_2$.  The perfect resemblance
between sums of spheres and the decomposition of the tensor
product is carried over to the multiplication of conjugacy
classes.  In fact, we set (note slight change of notation)
\begin{eqnarray*}
{\mathcal C}_{n}={\rm exp} (\Sigma_{\frac{n}{4\pi}}).
\end{eqnarray*}
Let $\frac{n}{4\pi}\equiv \alpha$ mod $\pi$ and
$\frac{m}{4\pi}\equiv \beta$ mod $\pi$.  By Theorem
\ref{thm:shaffaf1}
\begin{eqnarray*}
{\mathcal C}_n .{\mathcal C}_m= {\rm exp}
(\Sigma_{\frac{n}{4\pi}}+\Sigma_{\frac{m}{4\pi}})
\end{eqnarray*}
The representations corresponding to conjugacy classes contained
in ${\mathcal C}_n .{\mathcal C}_m$ correspond to those integers
$k\in [|m-n|,m+n]$ which $\equiv ~m+n$ mod 2.  Thus in the range
$[|m-n|,m+n]$ about half the representation occurring in
${\mathcal C}_n .{\mathcal C}_m$ appear in the Clebsch-Gordon
formula.  This is what is meant by the quantization of the product
of two conjugacy classes.

By the discretization of products of conjugacy classes ${\mathcal
C}_\alpha$ and ${\mathcal C}_\beta$ one means the choice of $N$
points on each and the determination of the corresponding
empirical measure of products of these points. If these $N$ points
are chosen randomly according to the invariant measures on the
conjugacy classes then the empirical measure of the products
converges weakly to the density function (\ref{eq:shaffaf1}) by
Theorem \ref{thm:shaffaf1} and is numerically demonstrated in
Figures 1 for a typical choice with $Np=2172$ points.

Each conjugacy class ${\mathcal C}_\alpha$ is naturally equivalent
to a copy of $S^2$.  It is therefore reasonable to investigate the
weak convergence of the empirical measure if the points on $S^2$
are chosen according to the requirements of Thomson's Problem of
distributing points on the sphere (see [KS]).  This means that the
points should be distributed so that the Coulomb potential
\begin{eqnarray*}
\sum_{i<j}^N \frac{1}{|z_i-z_j|^\alpha}
\end{eqnarray*}
is minimized.  A variation of this problem for $\alpha=1$ was
originally posed by J. J. Thomson in connection with his
investigations of the structure of the atom in 1904.  It remains
unsolved except for a few small values of $N$, and it has also
attracted attention for applications to complexity theory [Sm].
The lattice point method makes use of the natural embedding of the
icosahedron in $S^2$ and distributes $N=10(mn+m^2+n^2)+2$ points
on the sphere in such a way that the distribution exhibits a high
degree of symmetry and the points appear to be ``evenly" spaced.
It was conjectured in [A] that this distribution will provide the
solution to Thomson's Problem for $\alpha =1$.  The polar
coordinates method was devised in [KS] to achieve the minimum
required by Thomson's Problem and numerical tests disproved the
conjecture in [A] by showing that the (local) minimum achieved by
the polar coordinates method (where there was symmetry breakdown)
was in fact smaller.  It is therefore natural to test the
convergence of the empirical measure of products if the
discretization is done according the polar coordinates or the
lattice point methods.  In Figures 1 and 2 the convergence of the
empirical measure to the density predicted by Theorem
\ref{thm:shaffaf1} is exhibited for $Np=2172$ points.  The angles
$\phi$ and $\psi$ in the captions refer to the conjugacy classes
${\mathcal C}_\phi$ and ${\mathcal C}_\psi$ respectively. In
Figures 2 the corresponding empirical measures are calculated for
the lattice point method and it is noticed that even if the
measure converges to the density given by (\ref{eq:shaffaf1}), the
convergence is significantly slower. Similar conclusion is
applicable to the polar coordinates method as shown in Figures 3.

\begin{center}
\includegraphics[height=100mm]{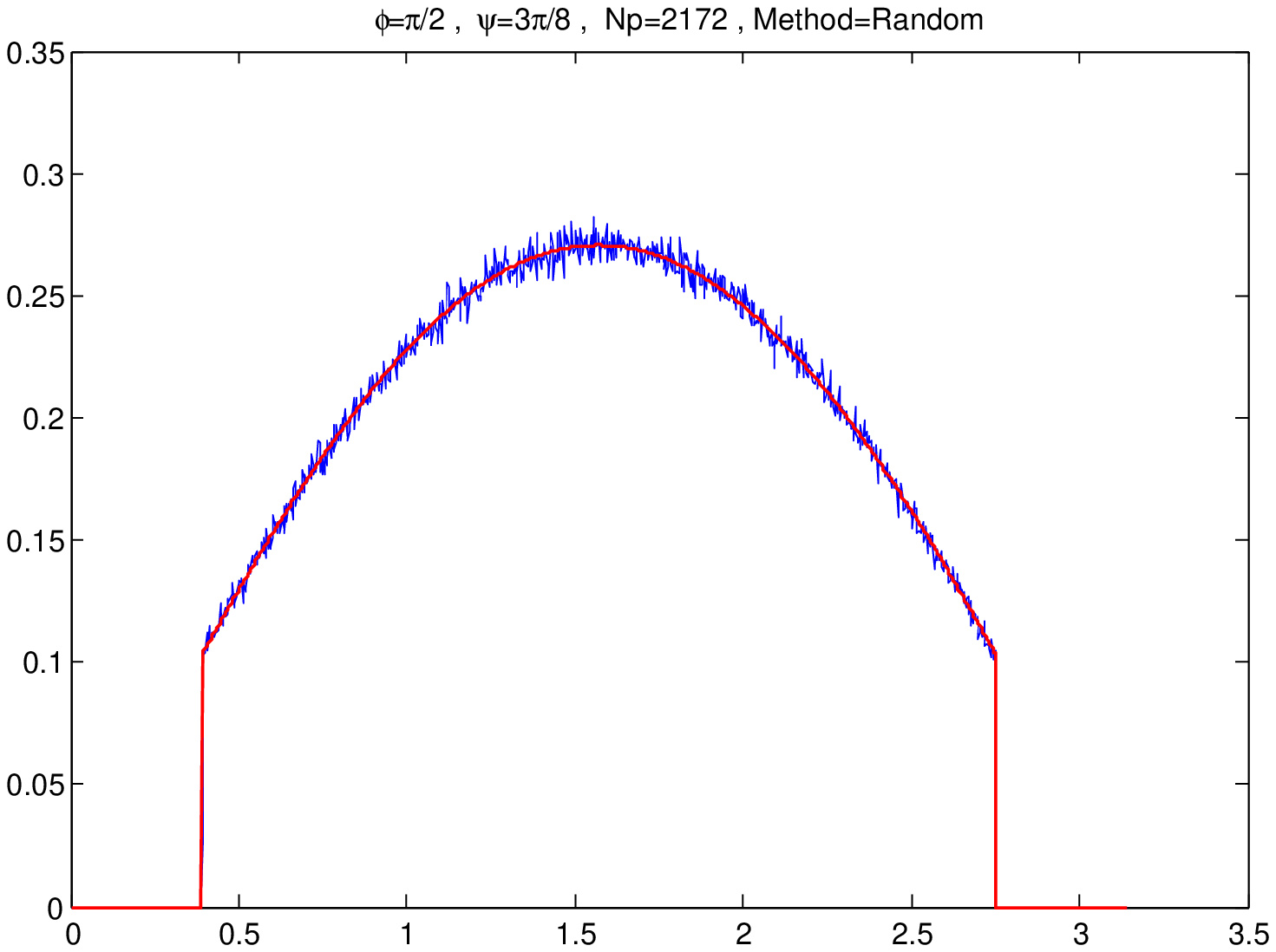} {\tiny Figure 1.}
\end{center}

\begin{center}
\includegraphics[height=100mm]{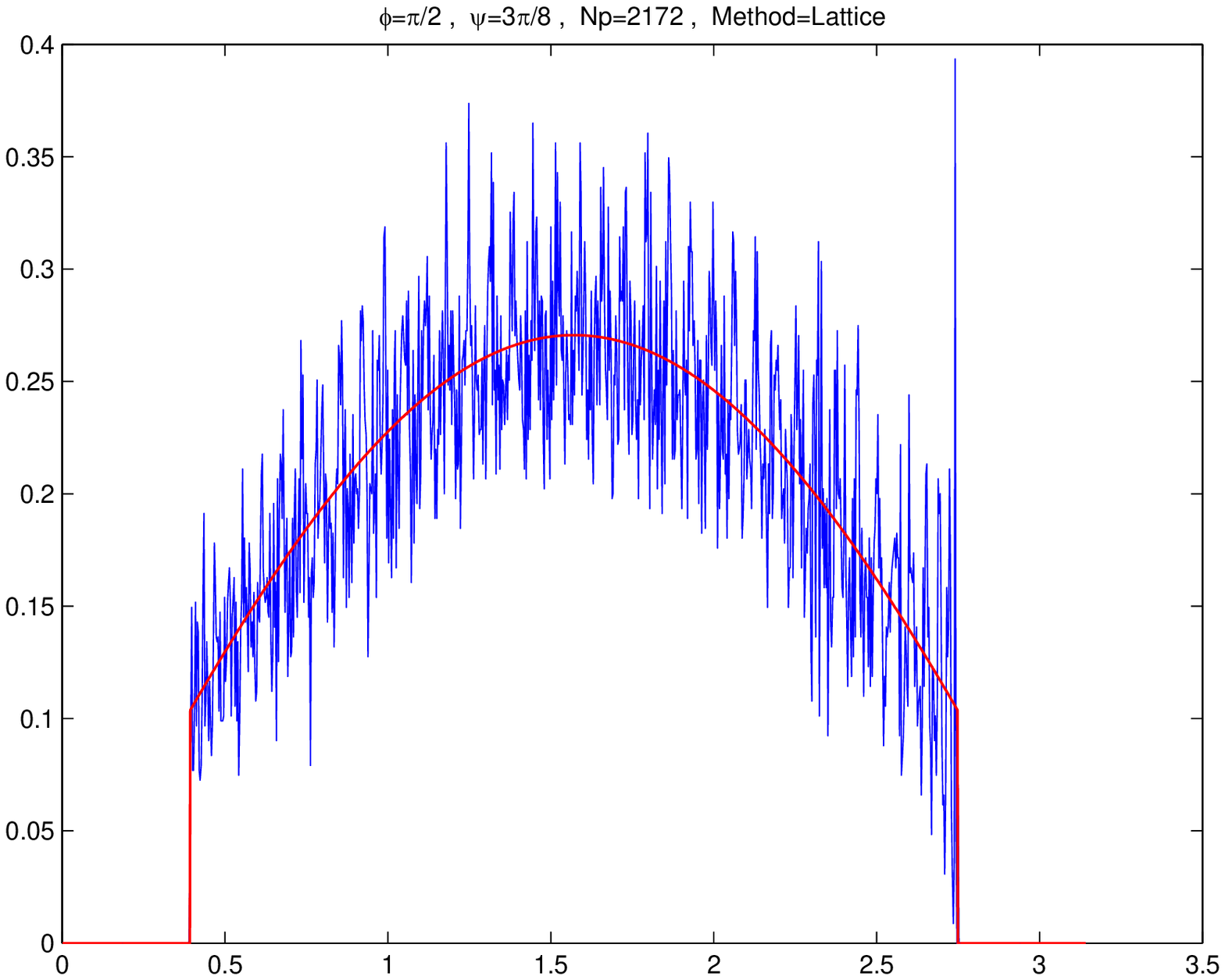} {\tiny Figure 2.}
\end{center}

\begin{center}
\includegraphics[height=85mm]{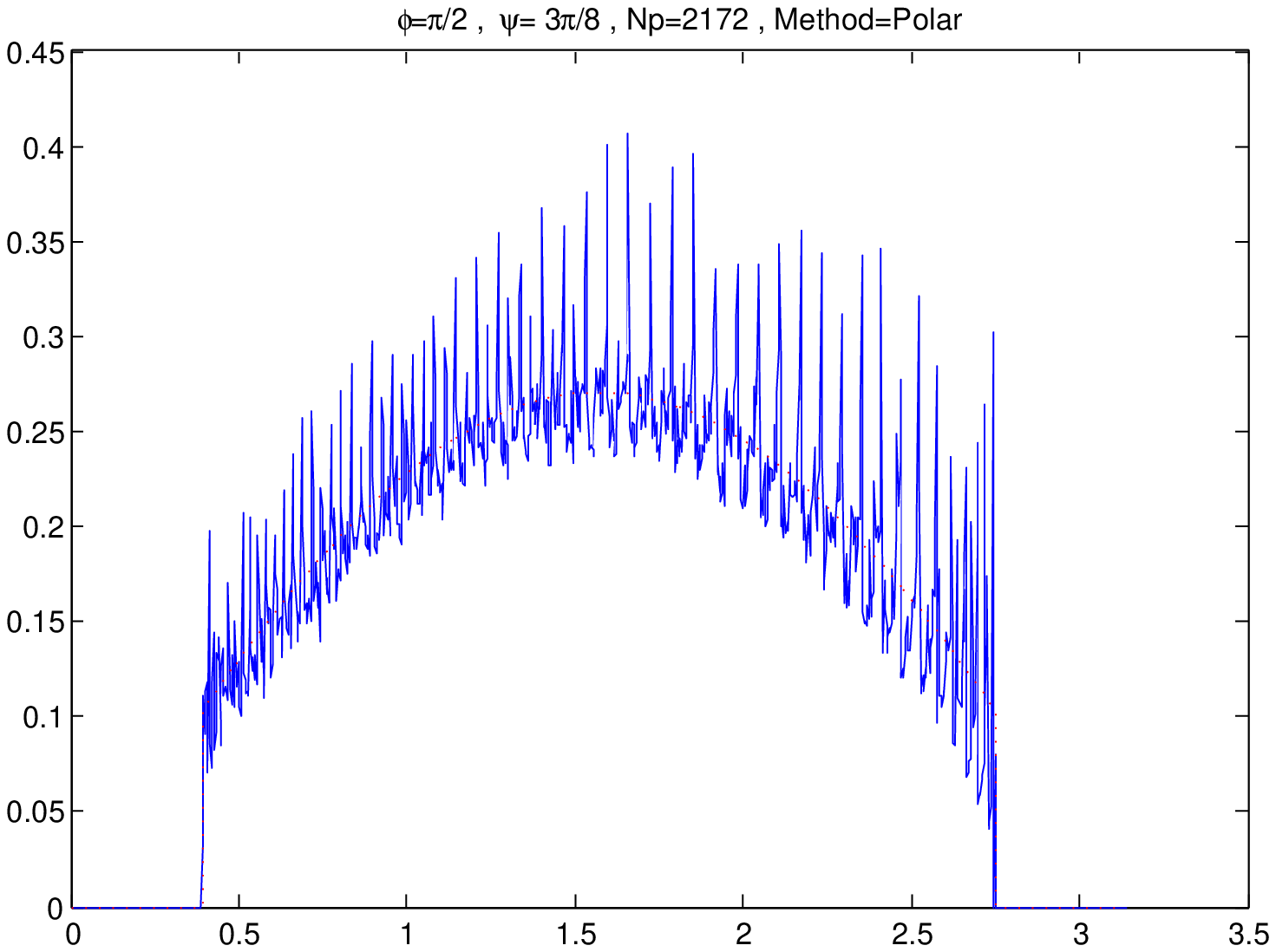} {\tiny Figure 3.}
\end{center}

\vfil\eject

\vskip .9 cm

\noindent Institute for Studies in Theoretical Physics and
Mathematics, Tehran, Iran, {\it and} \\ Sharif University of
Technology, Tehran, Iran. \vskip .4 cm

\end{document}